\documentstyle[11pt]{article}
\newcommand{\qed}{\hfill\rule{4pt}{8pt}\par\vspace{\baselineskip}}
\newtheorem{de}{Definition}[section]
\newtheorem{lm}[de]{Lemma}
\newtheorem{pr}[de]{Proposition}

\newtheorem{te}[de]{Theorem}

\newcommand{\otravez}{{1\le i \le \theta}}
\def\toba{{{\mathcal B}}}

\newcommand{\ydgg}{^{G}_{G}{\mathcal YD}}
\newcommand{\ydg}{^{k[G]}_{k[G]}{\mathcal YD}}
\newcommand{\yd}{^{\Gamma}_{\Gamma}{\mathcal YD}}
\newcommand{\ydga}{^{G(A)}_{G(A)}{\mathcal YD}}

\DeclareMathAlphabet{\mathfrak}{U}{euf}{m}{n}

\begin{document}


\title{Irreducible representations of liftings of quantum planes}
\author {N. Andruskiewitsch \thanks{This work was partially supported by CONICET,
Agencia C\'ordoba Ciencia, ANPCyT  and Secyt (UNC).}\\
Facultad de Matem\'atica, Astronom\'\i a y F\'\i sica,\\
 Universidad Nacional de C\'ordoba\\
CIEM -- CONICET \\
(5000) Ciudad Universitaria, C\'ordoba, Argentina\\  \small   email: andrus@mate.uncor.edu \\
 and \\ M.
Beattie \thanks{Research supported by NSERC.\newline 1991 {\it
Mathematics Subject Classification.}
Primary: 17B37. Secondary: 16W30} \\    Department of Mathematics and Computer Science \\
  Mount Allison University \\
   Sackville, NB, Canada E4L 1E6 \\ \small  email:
mbeattie@mta.ca}

\date{}\maketitle


\maketitle

\abstract{In this note, the irreducible representations of a
lifting of a quantum plane are determined. This paper appeared in:
 Lie theory and its applications in physics V, 414--423, World
 Sci. Publ., River Edge, NJ, 2004. }

\section{Introduction}

In \cite{beattieduals, ab}, the structure of the coradicals of the
duals of liftings of some quantum linear spaces was studied, and
several examples were explicitly constructed.  In this note, we
describe the irreducible representations for any lifting of a
quantum plane.  If $A$ is the lifting of a quantum linear plane,
then the approach in \cite{beattieduals, ab} was via the coradical
of $A^*$. The coalgebra $A^*=C$ was written as a direct sum of
certain sub-coalgebras $C(\eta)$ with $\eta \in \hat{\Gamma}$ and
the coradicals of the $C(\eta)$ were determined.  Here we write
the lifting $A$ as a product of algebras $A(\xi)$, $\xi \in
\hat{\Gamma}$, and determine irreducible representations of each
of six possible types of algebra $A(\xi)$ that can arise.

\subsection*{Notation}
Let $k$ be an algebraically closed field of characteristic 0. If
$G$ is a finite  group then $\ydgg = \ydg$ will denote the
category of Yetter-Drinfel'd modules over the group algebra
$k[G]$.

Throughout, we let  $\Gamma$ be a finite abelian group and let
$\hat{\Gamma}$ denote the group of  characters of $\Gamma$. For $V
\in \yd$, $g \in \Gamma, \chi \in \hat{\Gamma}$,   we write
$V^\chi_g$ for the set of $v \in V$ with the action of $\Gamma$ on
$v$ given by $h \to v =\chi(h)v$ and the coaction by $\delta (v) =
g \otimes v$. Since $\Gamma$ is an abelian group, $V= \oplus_{g\in
\Gamma, \chi \in \widehat{\Gamma}} V^\chi_g$ \cite[Section
2]{assurvey}.

The coradical of a coalgebra $C$ is denoted $C_0$.

If $v$ is a complex number and $a$ is a non-negative integer, then
we set as usual
$$
(0)_v = 0, \qquad (a)_v = 1 + v + \dots v^{a-1} = 1 + v(a -1)_v =
(a -1)_v + v^{a-1}.
$$

\section{Preliminaries}

Let $A$ be a finite-dimensional pointed Hopf algebra; let
$V$ be the infinitesimal braiding of $A$ \cite{assurvey}.
Recall that $V \in \ydga$ is the Yetter-Drinfeld module of right
coinvariant elements in the $k G(A)$-Hopf module $A_1/A_0$.

Let $\Lambda$ be a subgroup of $G(A)$ central in $A$. For $\xi \in
\widehat{\Lambda}$, let $$e_{\xi} = \vert \Lambda\vert ^{-1}
\sum_{g\in \Lambda} \xi^{-1}(g) g$$ denote the minimal idempotent
corresponding to $\xi$. If also $\eta \in \widehat{\Lambda}$
then
\begin{equation}\label{idemp}
e_{\eta}\leftharpoonup \xi
= \langle \xi,e_{\eta, (1)}\rangle e_{\eta, (2)} = e_{\xi^{-1}\eta}.
\end{equation}

\bigbreak Clearly, $Ae_{\xi}$ is a two-sided ideal of $A$ and $A =
\oplus _{\xi \in \widehat{\Lambda}} Ae_{\xi}$. Furthermore, let
$J_{\xi} = \oplus _{\xi\neq \zeta \in \widehat{\Lambda}}
Ae_{\zeta}$ and let $A(\xi) : = A/ J_{\xi}$; we have an
isomorphism of algebras $A \simeq \prod_{\xi \in
\widehat{\Lambda}} A(\xi)$. Hence, we have an isomorphism of
coalgebras $C := A^* \simeq \oplus_{\xi \in \widehat{\Lambda}}
C(\xi)$ where $C(\xi) := A(\xi)^* =J_{\xi}^{\bot} $. Clearly,
$C(\epsilon)$ is a Hopf subalgebra of $C$; it is dual to the
quotient Hopf algebra $A/ J_{\epsilon} = A/ A k[\Lambda]^+$. It
turns out that a useful approach to the description of the
coalgebra structure of $C$ is via the coalgebras $C(\xi)$ for a
well-chosen subgroup $\Lambda$. Indeed,  $C_0 :=  \oplus_{\xi \in
\widehat{\Lambda}} C(\xi)_0$ by \cite[9.0.1]{sw}.

\bigbreak Assume that $\xi\in \widehat{\Lambda}$ is the restriction
of $\widetilde \xi \in G(A)$. Then $\widetilde \xi C(\zeta) =  C
(\xi\zeta)$ for any $\zeta \in \widehat{\Lambda}$. Indeed,
$\widetilde \xi C(\zeta) \subset (Ae_{\eta})^{\bot} $ for any $\eta
\in \widehat{\Lambda}$, $\eta \neq \xi\zeta$, by (\ref{idemp}); thus
$\widetilde \xi C(\zeta) = J_{\xi\zeta}^{\bot}$.

\bigbreak Let us consider the following hypotheses.
\begin{enumerate}
\bigbreak
\item[(a)]  The group
$\Gamma := G(A)$ is abelian. Then there exists a basis $v_1, \dots
v_{\theta}$ of $V$ with $v_i \in V^{\chi_i}_{g_i}$ for all $i$.
Let $r_i > 1$ be the order of $\chi_i(g_i)$.

 \item[(b)] $V$ is a quantum linear space; that is,
$\chi_i (g_j) \chi_j (g_i) =1$ for $i \neq j$.

\end{enumerate}
\bigbreak Furthermore, assuming (a) and (b), we choose the
subgroup $\Lambda$ of $\Gamma$ as follows:
\begin{enumerate}
 \item[(c)]
 $\Lambda =\{ g \in \Gamma: \chi_i (g) =1$ for $1
\leq i \leq \theta \}.$
\end{enumerate}

 Note that by  the definition in (c), $\Lambda$ is central in $A$, and
 if $\chi_j^{r_j}= \epsilon$ or if $\chi_i\chi_j = \epsilon$ for
 all $i \neq j$, then
\begin{equation}\label{lambdacentral}
g_j^{r_j} \in \Lambda \mbox{ and  } \chi_j \in \widehat{\Gamma
/\Lambda}.
\end{equation}
   Also since the sequence $1
\to \widehat{\Gamma / \Lambda} \to \widehat{\Gamma} \to
\widehat{\Lambda} \to 1$ is exact, we can always choose  a
preimage $\widetilde \xi$ in $\widehat{\Gamma}$ of $\xi \in
\widehat{\Lambda}$; then $C(\xi) =\widetilde \xi C(\epsilon)$.

\bigbreak
\begin{pr}\label{liftqls}
Let $A$ be a finite-dimensional pointed Hopf algebra; let $V$ be
the infinitesimal braiding of $A$. Assume (a) and (b) above, i.e.
$A$ is a lifting of a quantum linear space with abelian group of
grouplikes. Then $A$ is generated by the grouplike elements $h\in
\Gamma$ and by $(1, g_i)$-primitives $x_i, 1 \leq i \leq \theta$
with defining relations
\begin{eqnarray*}
hx_i &=& \chi_i(h) x_ih; \\
x_i ^{r_i} &=& \alpha_{ii} (g_i^{r_i} -1);\\
x_i x_j &=& \chi_j (g_i) x_j x_i + \alpha_{ij} ( g_i g_j-1),
\end{eqnarray*}
for all $h\in \Gamma$, $\otravez$. We have dim $A = \vert \Gamma
\vert r_1 \dots r_\theta$.
\end{pr}
We may assume $\alpha_{ii} \in \{ 0,1\}$ and then we must have
that \begin{eqnarray*} \alpha_{ii} &=& 0 \mbox{ if } g_i^{r_i} =1
\mbox { or }\chi_i^{r_i} \neq \epsilon;\\
 \alpha_{ij} &=& 0 \mbox{ if }g_i g_j =1 \mbox{ or }\chi_i
\chi_j \neq \epsilon. \end{eqnarray*}
   Note that $\alpha_{ji} = - \chi_j (g_i)^{-1}
\alpha_{ij} = - \chi_i (g_j) \alpha_{ij}$. Thus the lifting $A$ is
described by the \emph{lifting matrix}
${\mathcal A} = ( \alpha_{ij})$ with 0's or 1's
on the diagonal and with $\alpha_{ji} =- \chi_i (g_j)
\alpha_{ij}$ for $i \neq j$.

\vspace{.25cm} {\bf Proof.} By the same argument as in \cite{asens},
$A$ is generated by group-like and skew-primitive elements. See
\cite{asp3} for the rest of the proof. These  liftings were
independently constructed in \cite{bdg} by iterated Ore extensions.
\qed \bigbreak

 In the rest of the paper we assume that
$A$ is a lifting of a quantum linear space with notation as in
Proposition \ref{liftqls}.

\bigbreak Since $\Gamma$ is finite abelian, there exist elements
$h_1, \dots, h_t$ in $\Gamma$, and non-negative integers $a_1,
\dots, a_t$,   with $a_u = $ ord $h_u$, $\Gamma = \langle h_1\rangle
\oplus \dots \oplus \langle h_t\rangle$.

Let $\xi \in \widehat{\Lambda}$; the algebra $A(\xi)$ is then
generated by $x_i$, $\otravez$, $h_u$, $1\le u \le t$, and relations
$$h_uh_{\ell} = h_{\ell}h_u,  \quad
h_ux_i = \chi_i(h_u) x_ih_u, \quad g = \xi(g) \mbox{ for } g \in
\Lambda,$$
$$x_i^{r_i} = \alpha_{ii}( g_i^{r_i}-1), \quad x_ix_j = \chi_j(g_i)
x_jx_i + \alpha_{ij} ( g_ig_j-1),
$$
$1\le u, \ell \le t$, $\otravez$. Let $V$ be any
finite-dimensional $A(\xi)$-module; we denote the action of the
elements $x_i$, $h_u$ on $V$ by the same letters. Thus $V$ is a
$\Gamma$-module, since $A(\xi)$ is a quotient of $A$; and $$V =
\oplus_{\eta \in F(\xi)} V^{\eta}$$ where
\begin{equation}\label{F}F(\xi)= \{\eta \in \widehat \Gamma: \eta
\vert_\Lambda = \xi \} = \widetilde  \xi \, \widehat{\Gamma /
\Lambda};\end{equation} see above. We have dim $A(\xi) = \frac{ |
\Gamma | r_1 \dots r_\theta}{| \Lambda |}$.

\bigbreak The case when the rank $\theta = 1$ is known \cite{rad,
beattieduals}. We shall investigate the case $\theta = 2$ in the
next Section. Consider the condition

\bigbreak\begin{enumerate}
\item[(d)] The rank $\theta > 2$ and ord $r_i > 2$, $\otravez$.
\end{enumerate}

\bigbreak
We say that $i, j \in I := \{1, \dots, \theta\}$ are \emph{linked}
\cite{asens} if $\alpha_{ij} \neq 0$. By (d), if $i$ is linked to
$j$ and $k$ then $j=k$ \cite{asens}. Thus $I$ is a disjoint union
of the set of vertices which are linked, which has even cardinal,
and the rest. Roughly speaking, the representation theory
of $A(\xi)$ looks like the ``tensor product'' of
representation theories of similar algebras with rank $\theta = 1$ or 2.

\bigbreak

\section{The rank 2 case}

In this section we assume that $\theta = 2$ and write $x = x_1$,
$y = x_2$, $r = r_1$, $s = r_2$. The lifting matrix of $A$ has the
form ${\mathcal A} = \left [
\begin{array}{cr} \alpha_{11} &\nu\\ -\chi_2(g_1)  \nu & \alpha_{22}
\end{array} \right]$.

Let $\xi \in \widehat{\Lambda}$; we shall determine the
irreducible representations of the algebra $A(\xi)$ generated by
$x$, $y$, $h_u$, $1\le u \le t$, and relations
\begin{eqnarray}
\label{relscero} h_uh_{\ell} &=& h_{\ell}h_u, \quad g = \xi(g)
\mbox{ for } g \in \Lambda;
\\ \label{relsuno} h_u x &=& \chi_1(h_u) x h_u;
\\ \label{relsdos} h_u y &=& \chi_2(h_u) y h_u ;
\\ \label{relstres} x^r &= & \alpha_{11}( g_1^r-1)  ;
\mbox{ we denote } \alpha_{11}(\xi(g_1^r)-1) = \alpha;
\\ \label{relscuatro} y^s &=& \alpha_{22}( g_2^s -1) ;
\mbox{ we denote } \alpha_{22}(\xi(g_2^s)-1)= \beta;
\\ \label{relscinco} xy &=& \chi_2(g_1) yx + \nu(g_1g_2-1).
\end{eqnarray}
We have dim $A(\xi) = \frac{\vert \Gamma \vert rs}{\vert \Lambda
\vert}$.

Let $q = \chi_1(g_1)$.  If $\chi_1\chi_2 = \epsilon$, then $q =
\chi_2(g_1)^{-1} = \chi_1(g_2) = \chi_2(g_2) ^{-1}$ (hence $r=s$),
and relation (\ref{relscinco}) becomes $yx = q\left(xy  -
\nu(g_1g_2-1) \right)$.

\bigbreak
We distinguish six cases; up to change of variables,
these six cases cover all possibilities for $A(\xi)$.

 \begin{enumerate}
 \item[\bf (I)]  $  \alpha = \beta = 0$, $\nu = 0$. \bigbreak
 \item[\bf (II)] $  \alpha = 1$, $\beta = 0$, $\nu = 0$.
Necessarily, $\chi_1^r = \epsilon$ and $g_1^r \neq 1$.

\bigbreak \item[\bf (III)] $\alpha = \beta = 1$, $\nu = 0$.
Necessarily, $\chi_1^r = \chi_2^s = \epsilon$, $g_1^r \neq 1$ and
$g_2^s \neq 1$.

\bigbreak \item[\bf (IV)] $\alpha = \beta = 0$, $\nu  = 1$.
Necessarily, $\chi_1\chi_2 = \epsilon$ and $g_1g_2 \neq 1$.

\bigbreak \item[\bf (V)] $\alpha = 1$, $\beta = 0$, $\nu = 1$.
Necessarily, $\chi_1^r = \epsilon$ and $g_1^r \neq 1$;
$\chi_1\chi_2 = \epsilon$ and $g_1g_2 \neq 1$.

\bigbreak \item[\bf (VI)] $\alpha = \beta = 1$, $\nu \neq 0$.
Necessarily, $\chi_1^r = \chi_2^s = \epsilon$, $g_1^r \neq 1$,
$g_2^s \neq 1$; $\chi_1\chi_2 = \epsilon$ and $g_1g_2 \neq 1$.
\end{enumerate}

Now we examine each of the 6 cases listed above.

\bigbreak {\bf Case (I)}  Here the lifting is trivial. In this
case we have the following well-known theorem.

\begin{te}The irreducible representations of $A(\xi)$ have dimension
one and are parametrized by $\widehat{\Gamma/\Lambda}$. Indeed,
$A(\xi) \simeq \toba(V) \# k[\Gamma] / (\toba(V) \#
k[\Gamma])e_{\xi}$.  \end{te}

\vspace{.5in}

{\bf Case (II)} Here $\alpha =1$ and $\beta = \nu = 0$.
 Thus $\chi_1^r = \epsilon$ and since
$q=\chi_1(g_1) $ is a primitive $r-$th root of unity, $r$ is the
order of $\chi_1$. \par It is now convenient to consider the
algebra $B(\xi)$ presented by generators $x$, $h_u$, $1\le u \le
t$, and relations
\begin{eqnarray*}h_uh_{\ell} = h_{\ell}h_u,   &\quad&  g = \xi(g) \mbox{ for } g \in \Lambda,
 \\ h_u x = \chi_1(h_u) x h_u,
&\quad & x^r = 1.
\end{eqnarray*}
The representation theory of $B(\xi)$   is well-known; we include
it for completeness. (For example, see \cite{yang}, \cite{rad}.)

\begin{lm}\label{basic}Let  $\eta\in F(\xi)$.
Let $W(\eta)$ be a vector space with a basis $f_i$, $0 \leq i \leq
r-1$, with subscripts $i \in {\mathbf Z}  /r$. There is a
representation of $B(\xi)$ on $W(\eta)$ defined by the following
rules:
$$ x.f_i = f_{i + 1}, \quad h_u.f_i = (\chi_1^i\eta)(h_u) f_i,
\quad 1\le u \le t$$ where $\quad i \in {\mathbf Z} /r$ means that
$x.f_{r-1}=f_0$. Furthermore, $W(\eta)$ is irreducible; all
irreducibles are of this kind; $W(\eta) \simeq W(\eta')$ if and only
if  $\eta = \eta' \chi^m_1$ for some $m$; and $B(\xi)$ is
semisimple.
\end{lm}

{\bf Proof.}  It is straightforward to verify that $W(\eta)$ is a
representation of $B(\xi)$. We see that $W(\eta)$ is irreducible
because the $f_i$'s belong to different isotypical components for
$\Gamma$.
\par Let $V$ be an irreducible representation of $A(\xi)$. Since
$V = \oplus_{\eta \in F(\xi)} V^{\eta}$, we can choose $v \in
V^{\eta} -0$ for some $\eta \in \widehat \Gamma$. Let $t$ be the
order of $x$ in $V$; clearly $x^i.v \in V^{\chi^i_1\eta} -0$,
$0\le i \le t$; and $x^t.v = v$. Then $r \vert t \vert r$,
thus $t=r$ and $V \simeq W(\eta)$. If $\phi: W(\eta) \to W(\eta')$
is an isomorphism of $B(\xi)$-modules, then $\phi(f_0) \in kf'_m$;
thus $\eta = \eta' \chi^m_1$. Conversely, assume that $\eta = \eta'
\chi^m_1$ and define $\phi: W(\eta) \to W(\eta')$ by $\phi(f_i) =
f'_{m + i}$; $\phi$ is an isomorphism of $B(\xi)$-modules. Finally,
$\dim B(\xi) \le \vert \Gamma / \Lambda\vert r$ by the defining
relations, but the dimension of the quotient of $B(\xi)$ by its
Jacobson radical is $\ge \vert \Gamma / \Lambda\vert r$ by what we
have already proved; thus $B(\xi)$ is semisimple. \qed

\bigbreak \begin{te}    Let $A(\xi)$ be such that $\alpha = 1$,
$\beta = \nu = 0$.  Let $\pi: A(\xi) \to B(\xi)$ be the algebra map
sending $y$ to $0$. Then any irreducible representation of $A(\xi)$
factorizes through $\pi$. In particular, all irreducible
representations of $A(\xi)$ are described by Lemma \ref{basic}.
\end{te}
\bigbreak {\bf Proof.} Let $V$ be an irreducible representation of
$A(\xi)$. Since $0 \neq \ker y$ is stable under the action of
$\Gamma$ and $x$, we see that $y$ acts as 0 on $V$. \qed

\vspace{.5in} {\bf Case (III)} Here $\alpha = \beta = 1$ and $\nu =
0$.
  Let  $w := \chi_2(g_1)$.  Since $\alpha =
  \beta = 1$, as noted in previous cases, we have that $\chi_1 $
  has order $r$ and $\chi_2$ has order $s$.  Thus $w^r = w^s = 1$.

 \begin{te} \label{caseiii} $A(\xi)$ is semisimple. \end{te}

\bigbreak {\bf Proof.} There exists a unique algebra automorphism
$Y$ of $B(\xi)$ such that\newline $Y(h_u) = \chi_2(h_u)^{-1} h_u$,
$1\le u \le t$, and $Y(x) = w^{-1}x$; clearly $Y^s = $id. It is
well-known that the smash product $B(\xi) \#  k\,Y$ is semisimple,
see \cite{MoW}. But $B(\xi) \#  k\,Y$ is isomorphic as an algebra
to $A(\xi)$, say by dimension counting.   \qed \smallbreak The
explicit description of all simple $A(\xi)$-modules can be
obtained by means of Clifford theory, see \cite{MoW}.

\vspace{.25in}
 {\bf Case (IV)}
In the next three cases, we have $\nu \neq 0$ so that $\chi_1
\chi_2 = \epsilon$ and so $r=s$.
 For all of the remaining cases, we will want to
define a set of scalars $c_i$ by the following recursive
definition.

Assume that $\chi_1\chi_2 = \epsilon $.   Fix $c = c_0$ and for
$\eta \in F(\xi)$ and $i> 0$, define
\begin{equation}\label{recursion}
c_i = q \left(c_{i-1} - \nu \eta \chi_1^{i-1}(g_1g_2 -1) \right)
 =  q \left(c_{i-1} + \nu - \nu q^{2(i-1)} \eta(g_1g_2) \right).
\end{equation}
The second equality follows from the fact that $q = \chi_1(g_1) =
\chi_1(g_2)$.    A simple induction shows that  for $i > 0$, we
have
\begin{equation}\label{closed} c_i = q^ic + q (i)_{q}\nu \left( 1- q^{i-1}
\eta(g_1g_2)\right).\end{equation}

Thus if $q$ is a primitive $r$-th root of unity, if $i \equiv k
\mbox{ mod } r$, then $c_i = c_k$.

\vspace{.2in} In this case, $\alpha = \beta = 0$ and $\nu = 1$.
 Here, the representation theory is similar to that of a
Frobenius-Lusztig kernel of type $  \mathfrak{sl}$(2).

\bigbreak \begin{te}Let $\eta \in F(\xi)$ as defined in (\ref{F}).
Let
  $(c_i)_{i\ge 0}$
be  scalars defined recursively by (\ref{recursion}) with $c_0=0$.

Let $N$ be the least positive integer such that $c_N =0$. Note
that since $(r)_q=0$, then $N \leq r$. Let $L(\eta)$ be a vector
space with a basis $(v_i)_{0\le i \le N-1}$; set $ v_{-1} = v_{N}
= 0$ in $L(\eta)$. Then there exists a representation of $A(\xi)$
on $L(\eta)$ given by
\begin{equation}\label{IIdos}
h_u.v_i = \eta\chi_1^{i}(h_u) v_i, \quad 1\le u \le t, \qquad
y.v_i = c_i v_{i-1}, \qquad x.v_{i} =  v_{i + 1}.
\end{equation}
Furthermore, $L(\eta)$ is irreducible.   Also any irreducible
$A(\xi)$-module is isomorphic to $L(\eta)$ for some $\eta$; and
$L(\eta)$ is isomorphic to $L(\eta')$ only when $\eta = \eta'$.
\end{te}

\bigbreak {\bf Proof.}  The verification that (\ref{IIdos}) defines
a representation of $A(\xi)$ is straightforward. The fact that $
y.v_{0} = 0 = x.v_{N-1} $  ensures that relations (\ref{relstres})
and (\ref{relscuatro}) hold while (\ref{recursion}) guarantees that
relation (\ref{relscinco}) is respected.  We leave the reader to
check the details and to check that $A(\xi)$ is irreducible.
\par Let $\rho: A(\xi) \to$ End $V$ be an irreducible
representation. Since $\ker y \neq 0$ and is $\Gamma$-stable, there
exists $v\in \ker y -0$, $v\in V^{\eta}$ for some $\eta \in F(\xi)$.
Set $v_0 = v$, $v_i = x^i.v $, $i> 0$; $v_i\in V^{\eta\chi_1^{i}}$
for all $i$. It follows from the fact that relation
(\ref{relscinco}) must be respected and from a simple induction that
for $i > 0$, we have $y.v_i = d_i v_{i-1}$, where the $d_i$'s
satisfy the recursive relation (\ref{recursion}).

Now, $v_0, v_1 \ldots v_{m-1}$ generate a submodule of $V$ where
$v_{m} = x^{m}.v_0 = 0$ and $m \leq r$.  If $m > N$, then
$v_{m-1}, \ldots , v_N$ is a submodule of $V$.  Thus $m=N$ and
since $d_i = c_i$, $i\ge 0$, then $V$ coincides with the submodule
generated by the $v_i$'s, $i\ge 0$ which is isomorphic to
$L(\eta)$.

Finally, $L(\eta)$ is presented as $A(\xi)$-module by generator
$v_0$ with relations $h_u(v_0) = \eta(h_u) v_0$, $1\le u \le t$,
$y.v_0 = 0$, $x^N.v_{0} = 0$. Thus, $L(\eta)\simeq L(\eta')$ implies
$\eta = \eta'$. \qed

\vspace{.25in} {\bf Case (V)} In this case, we have $\alpha = \nu
= 1$ and $\beta =0$.

 \begin{te} Let  $\eta\in F(\xi)$. Let
$W(\eta)$ be the $B(\xi)$-module defined in Lemma \ref{basic} with
basis $f_0, \ldots, f_{r-1}$ with subscripts taken modulo $r$.
   Set $c_i = 0 $ for $i = 0$ and
define scalars $( c_i)_{0 < i\le r}$ recursively  by
(\ref{recursion}).
Define an operator $y$ on  $W(\eta)$ by $y.f_i = c_if_{i-1}$. Then
this defines a representation of $A(\xi)$  and we denote this
$A(\xi)$-module by $L(\eta)$. Then $L(\eta)$ is irreducible; all
irreducibles are of this kind; If $L(\eta) \simeq L(\eta')$ then
$\eta = \eta' \chi^m_1$ for some $m$. \end{te}

\bigbreak {\bf Proof.} As usual, the verification that $L(\eta)$ is
an irreducible representation is straightforward. Recall from
(\ref{closed})  that it makes sense to compute subscripts modulo
$r$.


Next, let $V$ be an irreducible $A(\xi)$-module. Then there exists
$\eta\in F(\xi)$ and $v\neq 0$ such that $v\in \ker y \cap
V^{\eta}$. Set $f_0 = v$, $f_{i} = x^iv$. Arguing as in Lemma
\ref{basic} we see that the $f_i$'s span a $B(\xi)$-submodule $U$
isomorphic to $W(\xi)$. Relation (\ref{relscinco}) implies the
description of the action of $y$ on $U$ by $y.f_i =
\alpha_if_{i-1}$, where the $\alpha_i$'s are defined by
(\ref{recursion}). Hence $V = U \simeq L(\xi)$. \qed

\vspace{.25cm} {\bf Case (VI)}  In this case, $\alpha = \beta = 1$
and $\nu \neq 0$ .

 \begin{te}Let  $\eta\in F(\xi)$. Let
$W(\eta)$ be the $B(\xi)$-module defined in Lemma \ref{basic}. Set
$c_0 = c \neq 0$   and define a family of scalars $(c_i)_{0 \leq i
\leq r-1}$ inductively by (\ref{recursion}). Define an operator
$y$ on  $W(\eta)$ by $y.f_i = c_if_{i-1}$.

(i). This defines a representation of $A(\xi)$ if and only if $c$
is a solution of the equation
\begin{equation}
\label{VIdos}
 c_0c_{1}\dots c_{r - 1} = 1.
\end{equation}
If this is the case, we denote this $A(\xi)$-module by $L(\eta,
c)$. Furthermore, $L(\eta, c)$ is irreducible.

(ii).  $L(\eta, c) \simeq L(\eta', c')$ if and only if $\eta = \eta'
\chi^m_1$ for some $m$, and $c = c'_m$.

(iii). Assume that the equation (\ref{VIdos}), a polynomial in $c$,
has simple roots. Then all irreducibles are of this kind and
$A(\xi)$ is semisimple.
\end{te}

\bigbreak {\bf Proof.} (i). We have to check the relations
(\ref{relsdos}), (\ref{relscuatro}) and (\ref{relscinco}). Here
(\ref{relsdos}) is clear and (\ref{relscuatro})  is equivalent to
(\ref{VIdos}). We evaluate both sides of (\ref{relscinco}) on $f_i$;
if $0\le i < r-1$ the equality follows from the defining condition
(\ref{recursion}); otherwise it follows from (\ref{closed}). The
irreducibility is clear.

(ii). Left to the reader.

(iii). In general, the dimension of the semisimple quotient of
$A(\xi)$ corresponding to all the representations of the type
$L(\eta, c)$ is $\frac{| \Gamma | \#  \{\mbox{solutions of }
(\ref{VIdos})\} r^2}{\vert \Lambda \vert r}$, and this equals $\dim
A(\xi)$ if and only if the equation (\ref{VIdos}), a polynomial in
$c$, has simple roots.
 \qed


\section{CONCLUSIONS}

Let $A$ be a lifting of a quantum plane. Then $A \simeq \prod_{\xi
\in \widehat{\Lambda}} A(\xi)$; hence the irreducible
representations of $A$ are the union of the irreducible
representations of $A(\xi)$, $\xi \in \widehat{\Lambda}$. We have
determined the last ones in Section 3, up to a finite number of
exceptions in Case (VI). As a consequence, we can also determine the
coradical of the dual Hopf algebra $C = A^*$.


\section*{Acknowledgments}

The first author thanks Sonia Natale for conversations on Case
(III). He also thanks Vlado Dobrev, Toshko Popov and all the local
organizers of the V. International Workshop LIE THEORY AND ITS
APPLICATIONS IN PHYSICS for the beautiful organization. Both
authors thank Hans-J\"urgen Schneider for pointing out a mistake
in the published version of the paper.


\end{document}